\documentclass[12pt]{amsart}
\usepackage[latin1]{inputenc}
\usepackage{graphics}
\usepackage{graphicx}

\newtheorem{theorem}{Theorem}%[section]
\newtheorem{lemma}[theorem]{Lemma}
\newtheorem{proposition}[theorem]{Proposition}
\newtheorem{corollary}[theorem]{Corollary}

\theoremstyle{definition}

\theoremstyle{remark}

\usepackage{amssymb,amsmath}

\def\H{\mathbb{H}}

\newcommand{\m}{\mbox{$M$}}
\newcommand{\n}{\mbox{$\nabla$}}

\newcommand{\s}{\mbox{$\Sigma$}}
\newcommand{\R}{\mbox{${\mathbb R}$}}
\newcommand{\N}{\mbox{$\m^2\times\R_1$}}
\newcommand{\g}[2]{\mbox{$\langle #1 ,#2 \rangle$}}
\newcommand{\fle}{\mbox{$\rightarrow$}}
\newcommand{\rf}[1]{\mbox{(\ref{#1})}}
\newcommand{\rl}[1]{{~\ref{#1}}}

\newcommand{\fs}{\mbox{$\mathcal{C}^\infty(\s)$}}
\newcommand{\f}{\mbox{$f:\Sigma^2\fle\N$}}

\def\beq{\begin{equation}}
\def\eeq{\end{equation}}

\begin{document}

\title[Parabolicity of maximal surfaces]
{Parabolicity of maximal surfaces in Lorentzian product spaces}

\author{Alma L. Albujer}
\address{Departamento de Matem\'{a}ticas, Universidad de Murcia, E-30100 Espinardo, Murcia, Spain}
\email{albujer@um.es}

\author{Luis J. Al\'\i as}
\address{Departamento de Matem\'{a}ticas, Universidad de Murcia, E-30100 Espinardo, Murcia, Spain}
\email{ljalias@um.es}
\thanks{This work was partially supported by MEC projects MTM2007-64504 and Fundaci\'{o}n S\'{e}neca
project 04540/GERM/06, Spain.
This research is a result of the activity developed within the framework of the Programme in Support of Excellence Groups of the Regi\'{o}n de Murcia, Spain, by Fundaci\'{o}n S\'{e}neca, Regional Agency
for Science and Technology (Regional Plan for Science and Technology 2007-2010).}

%    General info

%\date{August 2008}%{First version, October 2007}

\dedicatory{Dedicated to Professor Marcos Dajczer on the occasion of his 60th birthday}


\begin{abstract}
In this paper we establish some parabolicity criteria for maximal surfaces immersed into a Lorentzian product space of the
form $\m^2\times\R_1$, where $\m^2$ is a connected Riemannian surface with non-negative Gaussian curvature and
$\m^2\times\R_1$ is endowed with the Lorentzian product metric $\g{}{}=\g{}{}_M-dt^2$. In particular, and as an
application of our main result, we deduce that every
maximal graph over a starlike domain $\Omega\subseteq M$ is parabolic. This allows us to give an alternative proof of the
non-parametric version of the Calabi-Bernstein result for entire maximal graphs in \N.
\end{abstract}

\maketitle

\section{Introduction}

A \textit{maximal} surface in a $3-$dimensional Lorentzian manifold is a spacelike surface with zero mean curvature.
Here by \textit{spacelike} we mean that the induced metric from the ambient Lorentzian metric is a Riemannian metric on
the surface. The mathematical interest of maximal surfaces resides in the fact that they locally maximize area among all
nearby surfaces having the same boundary \cite{Fr,BF}. Furthermore, maximal surfaces and, more generally spacelike
surfaces with constant mean curvature, have also a great importance in General Relativity \cite{MT}.

Parabolicity is a concept that lies in the borderline between several branches of mathematics, such as
Riemannian Geometry, Stochastic Analysis, Partial Differential Equations and Potential Theory.
Let us recall that a Riemannian surface
$(\Sigma,g)$ with non-empty boundary, $\partial\Sigma\neq\emptyset$, is said to be \textit{parabolic} if every bounded
harmonic function on $\Sigma$ is determined by its boundary values. It is interesting to observe that the parabolicity of
$\Sigma$ is equivalent to the existence of a proper nonnegative superharmonic function on $\Sigma$ (see beginning of Section\rl{s3}; for details see
\cite{MP} and \cite{P}).

Classically, a Riemannian surface without boundary is called parabolic if it does not admit
a nonconstant negative subharmonic function. In fact, in \cite{AA} we considered that
definition of parabolicity. However, along this work we will reserve the term \textit{parabolicity}
for Riemannian surfaces with non-empty boundary, and we will use the term \textit{recurrence} for
Riemannian surfaces without boundary. A Riemannian surface $(\Sigma,g)$ without boundary is
called \textit{recurrent} if for every nonempty open set $O\subset\Sigma$ with smooth boundary,
$\Sigma\setminus O$ is parabolic (as a surface with boundary). It is worth pointing out that the concept of recurrence
is equivalent to the classical concept of parabolicity (for surfaces without boundary). This follows from the observation that
a Riemannian manifold without boundary is recurrent precisely when almost all Brownian paths are dense in the manifold
(for the details, see \cite{MP} and \cite[Theorem 5.1]{Gr}).

In \cite{FL} Fern\'{a}ndez and L\'{o}pez have recently proved that properly immersed maximal surfaces with non-empty boundary in
the Lorentz-Minkowski spacetime $\R_1^3$ are parabolic if the Lorentzian norm on the maximal surface in $\R_1^3$ is
eventually positive and proper (see also Section 4 in \cite{FL2}). Motivated by that work, in this paper we study some parabolicity criteria for maximal
surfaces immersed into a Lorentzian product space of the form $\m^2 \times \R$, where $\m^2$ is a connected Riemannian
surface and $\m^2 \times \R$ is endowed with the Lorentzian metric
\[
\g{}{}=\pi^\ast_M(\g{}{}_M)-\pi^\ast_{\mathbb{R}}(dt^2).
\]
Here $\pi_M$ and $\pi_\mathbb{R}$ denote the projections from $M \times \R$ onto each factor, and $\g{}{}_M$ is the
Riemannian metric on \m. For simplicity, we will simply write
\[
\g{}{}=\g{}{}_M-dt^2,
\]
and we will denote by $\m^2 \times \R_1$ the 3-dimensional product manifold $\m^2 \times \R$ endowed with that Lorentzian
metric. Observe that in particular, when $\m^2=\R^2$ is the flat Euclidean plane what we obtain is the well known
$\N=\R_1^3$ Lorentz-Minkowski 3-dimensional spacetime.

A natural generalization of the Lorentzian norm on a surface in $\R_1^3$ to the Lo\-rent\-zian product \N\ consists in
considering the function $\phi=r^2-h^2$ where the function $r$ measures the distance on the factor $\m$ to a
fixed point $x_0 \in \m$ and $h \in \fs$ is the height function of the surface $\s$ (for the details, see Section\rl{s3}).
In this context, our main result is Theorem \ref{teorlap} which asserts that given a complete Riemannian surface $M^2$
with non-negative Gaussian curvature, then every maximal surface in \N\ with non-empty boundary and such that the
function $\phi:\s \fle \R$ is eventually positive and proper is parabolic. In particular, and as an application of this
result, we deduce that every maximal graph over a starlike domain $\Omega\subseteq M$ is parabolic. This
allows us to give an alternative proof of the non-parametric version of our Calabi-Bernstein result
\cite[Theorem 4.3 and Corollary 4.4]{AA} for maximal entire graphs in \N\ (see Corollary \ref{CB}).

The authors would like to heartily thank V. Palmer, J. P\'{e}rez and  S. Verpoort for several explanations and useful comments during the preparation
of this paper.

\section{Preliminaries}

A smooth immersion $f:\s^2 \rightarrow \m^2 \times \R_1$ of a connected surface $\s^2$ is
said to be a spacelike surface if $f$ induces a Riemannian metric on \s, which as usual is
also denoted by $\g{}{}$. In that case, since
\[
\partial_t=(\partial/\partial_t)_{(x,t)}, \quad x \in M, \, t \in \R,
\]
is a unitary timelike vector field globally defined on the ambient spacetime $\m^2 \times
\R_1$, then there exists a unique unitary timelike normal field $N$ globally defined on \s\
which is in the same time-orientation as $\partial_t$, so that
\[
\g{N}{\partial_t} \leq -1 < 0 \quad \textrm{on} \quad \s.
\]
We will refer to $N$ as the future-pointing Gauss map of \s, and we will denote by
$\Theta:\s \fle \left( -\infty,-1 \right]$ the smooth function on \s\ given by
$\Theta=\g{N}{\partial_t}$. Observe that the function $\Theta$ measures the hyperbolic
angle $\theta$ between the future-pointing vector fields $N$ and $\partial_t$ along \s.
Indeed, they are related by $\cosh \theta=-\Theta$.

Let $\overline{\n}$ and $\n$ denote the Levi-Civita connections in \N\ and \s,
respectively. Then the Gauss and Weingarten formulae for the spacelike surface \f\ are
given by
\beq
\label{eqgauss}
\overline{\n}_XY=\n_XY-\g{AX}{Y}N
\eeq
and
\beq
\label{eqwein}
AX=-\overline{\n}_XN,
\eeq
for every tangent vector fields $X,Y \in T\s$.
Here $A:T\s \fle T\s$ stands for the shape operator (or second fundamental form) of \s\
with respect to its future-pointing Gauss map $N$.

The height function $h \in \fs$ of a spacelike surface $f:\s^2 \fle \N$ is the projection
of \s\ onto $\R$, that is, $h=\pi_\mathbb{R} \circ f$. Observe that the gradient of
$\pi_\mathbb{R}$ on \N\ is
\[
\overline{\n} \pi_\mathbb{R}=-\g{\overline{\n} \pi_\mathbb{R}}{\partial_t}\partial_t=-\partial_t,
\]
so that the gradient of $h$ on $\Sigma$ is
\[
\n h=\left(\overline{\n} \pi_\mathbb{R}\right)^{\top}=-\partial_t^\top.
\]
Throughout this paper, for a given vector field $Z$ along the immersion, we will denote by $Z^\top \in T\s$ its
tangential component; that is,
\[
Z=Z^\top-\g{N}{Z}N.
\]
In particular, $\nabla h=-\partial_t-\Theta N$ and we easily get
\beq
\label{normgradh}
\|\n h\|^2=\Theta^2-1,
\eeq
where $\|\cdot\|$ denotes the norm of a
vector field on \s. Since $\partial_t$ is parallel on \N\ we have that
\beq
\label{eqaux}
\overline{\n}_X \partial_t=0
\eeq
for every tangent vector field $X \in T\s$. Writing
$\partial_t=-\nabla h-\Theta N$ along the surface \s\ and using Gauss \eqref{eqgauss} and Weingarten
\eqref{eqwein} formulae, we easily get from \eqref{eqaux} that
\beq
\label{eqaux2}
\n_X \n h=\Theta AX
\eeq
for every $X \in T\s$. Therefore the Laplacian on \s\ of the height
function is given by
\[
\Delta h=\Theta \mathrm{tr}A=-2H\Theta,
\]
where
$H=-(1/2) \mathrm{tr}A$ is the mean curvature of \s\ relative to $N$. In particular,
\beq
\label{laph}
\Delta h=0
\eeq
for every maximal surface in \N.

Given any function $\hat{\psi} \in \mathcal{C}^\infty(M)$, we can consider its lifting
$\bar{\psi}\in\mathcal{C}^\infty(\N)$ defined by
\[
\bar{\psi}(x,t)=\hat{\psi}(x).
\]
In turn, we associate to $\hat{\psi}\in\mathcal{C}^\infty(M)$ the function $\psi\in\mathcal{C}^\infty(\s)$ given by
$\psi=\bar{\psi} \circ f$. In this context, the Laplacian on \s\ of $\psi$ can be
expressed in terms of the Laplacian $\bar{\Delta}$ of $\bar{\psi}$ and the differential operators of $\hat{\psi}$ as follows.
\begin{lemma}
\label{lemmaaux}
Along a spacelike surface \f\ we have that
\beq
\label{lemma}
\Delta \psi=\bar{\Delta} \bar{\psi}+2H\g{N^\ast}{\hat{\n}\hat{\psi}}_M+\hat{\n}^2 \hat{\psi}(N^\ast,N^\ast)
\eeq
where $N^\ast=\pi_M^\ast(N)=N+\Theta\partial_t$, and $\hat{\n}$ and $\hat{\n}^2$ denote the gradient and the Hessian
operators on $M$, respectively.
\end{lemma}
\begin{proof}
Since $\bar{\n}\bar{\psi}=\n \psi-\g{\bar{\n}\bar{\psi}}{N}N$, we get from \eqref{eqgauss}
and \eqref{eqwein} that the Hessian operators of $\bar{\psi}$ and $\psi$ satisfy
\[
\bar{\n}^2\bar{\psi}(X,X)=\n^2\psi(X,X)+\g{AX}{X}\g{\bar{\n}\bar{\psi}}{N}
\]
for every $X \in T\s$. Therefore, it can be easily seen that
\beq
\label{rellap}
\bar{\Delta}\bar{\psi}=\Delta \psi-2H\g{\bar{\n}\bar{\psi}}{N}-\bar{\n}^2\bar{\psi}(N,N).
\eeq
Observe now that, as the function $\bar{\psi}$ does not depend on $t$, then
$\bar{\n}\bar{\psi}(x,t)=\hat{\n}\hat{\psi}(x)$. Thus, $\bar{\n}_N
\bar{\n}{\bar{\psi}}=\hat{\n}_{N^\ast}\hat{\n}{\hat{\psi}}$ and
\[
\bar{\n}^2 \bar{\psi} (N,N)=\hat{\n}^2\hat{\psi} (N^\ast,N^\ast),
\]
so that the lemma follows directly from \eqref{rellap}.
\end{proof}

\section{Parabolicity of maximal surfaces}
\label{s3}
Our main result in this section generalizes \cite[Theorem 3.1]{FL} (see also \cite[Theorem 4.2]{FL2}) to the case of maximal spacelike surfaces in
$M \times \R_1$, when \m\ is a complete Riemannian surface with non-negative Gaussian curvature. In that case, consider the
function $\hat{r}:M \fle \R$ defined by $\hat{r}(x)=\mathrm{dist}_M(x,x_0)$ where $x_0 \in M$ is a fixed point. Observe
that $\hat{r} \in \mathcal{C}^\infty(M)$ almost everywhere. Actually, $\hat{r}$ is smooth on $M\setminus\mathrm{Cut}(x_0)$, where
$\mathrm{Cut}(x_0)$ stands for the cut locus of $x_0$. As is well-known, $\mathrm{dim}\mathrm{Cut}(x_0)<2$ and $\mathrm{Cut}(x_0)$ is a null set.

Following our notation above, let $\bar{r}(x,t)=\hat{r}(x)$ denote the
lifting of $\hat{r}$ to \N, and for a given spacelike surface \f, let $r$ stands for the restriction of $\bar{r}$ to
$\Sigma$, $r=\bar{r}\circ f$. Consider $\Pi=\pi_M\circ f:\Sigma\fle M$ the projection of $\Sigma$ on $M$. It is not difficult
to see that $\Pi^*(\g{}{}_{M})\geq\g{}{}$, where $\g{}{}$ stands for the Riemannian metric on \s\ induced
from the Lorentzian ambient space. That means that $\Pi$ is a local diffeomorphism which increases the distance
between the Riemannian surfaces \s\ and \m, and by \cite[Chapter VIII, Lemma 8.1]{KN} $\Pi$ is a covering map(for the details, see Lemma 3.1 in \cite{AA}).
Therefore, $\mathrm{dim}\,\Pi^{-1}(\mathrm{Cut}(x_0))=\mathrm{dim}\mathrm{Cut}(x_0)<2$ and the function $r$ is smooth almost everywhere in \s.

The proof of our main result is based on the following technical lemma.
\begin{lemma}
\label{lemadebil}
Let $\Sigma^2$ be a Riemannian surface with non-empty boundary, $\partial \s \neq\emptyset$. If there exists a proper continuous function $\psi:\s\rightarrow\mathbb{R}$ which is eventually positive and superharmonic, then
$\Sigma$ is parabolic.
\end{lemma}
As is usual, by eventually we mean here a property that is satisfied outside a compact set.
\begin{proof}
The proof follows the ideas of the proof of an analogous criterium for proper smooth functions given by Meeks and P\'{e}rez \cite{MP,P}. For the sake of completeness, we sketch it here. It suffices to see that if $\varphi$ is a bounded
harmonic function on \s\ which vanishes on the boundary, $\varphi|_{\partial\Sigma}\equiv 0$, then $\varphi$ is constant zero. Let $K\subset\Sigma$ be a compact subset such that $\psi$ is positive and superharmonic on $\Sigma\setminus K$.
It suffices to show that $\varphi|_{\Sigma\setminus K}\equiv 0$. Let us assume that there exists a point $p_0\in\Sigma\setminus K$ such that $\varphi(p_0)\neq 0$. Since $\psi(p_0)>0$, there is a constant $a\in\mathbb{R}$ with
$a\varphi(p_0)>\psi(p_0)>0$. Since $\varphi$ is harmonic, the function  $\eta=\psi-a\varphi:\Sigma\setminus K\rightarrow\mathbb{R}$ is superharmonic on $\Sigma\setminus K$. Besides, since $\psi$ is proper and positive on
$\Sigma\setminus K$, then $\eta$ is positive outside the compact subset $K'=\{ p\in\Sigma\setminus K : \psi(p)\leq C\}$ where $C=\sup_{p\in\Sigma\setminus K}a\varphi(p)>0$. Finally, since $\eta(p_0)<0$, $\eta$ must reach its
minimum on $\Sigma\setminus K$ at an interior point, which is a contradiction by the minimum principle for superharmonic functions.
\end{proof}

Now we are ready to prove our main result.
\begin{theorem}
\label{teorlap}
Let $M^2$ be a complete Riemannian surface with non-negative Gaussian curvature.
Consider $\Sigma$ a maximal surface in \N\ with non-empty boundary, $\partial \s \neq\emptyset$, and assume that the function
$\phi:\Sigma \fle \R$ defined by
\[
\phi(p)=r^2(p)-h^2(p)
\]
is eventually positive and proper. Then $\Sigma$ is parabolic.
\end{theorem}

It is worth pointing out that
the assumption on the non negativity of the Gaussian curvature of $M$ is necessary. Actually, let $M^2=\H^2$ and consider
$\Omega\subset\H^2$ a connected domain with smooth boundary. Then, for a fixed $t_0\in\R$,
$\Sigma_{t_0}=\{(x,t_0)\in\H^2\times\R : x\in\Omega \}$ is trivially a non-parabolic maximal surface in
$\H^2\times\R$ on which $\phi$ is eventually positive and proper.

\begin{proof}
Let $a>1$, and consider $K=\{ p \in \s:\phi(p) \leq a \}\subseteq\s$. $K$ is compact because $\phi$
is eventually positive and proper. As is well known, parabolicity is not affected by adding or removing compact subsets,
so that \s\ is parabolic if and only if $\s\setminus K$ is parabolic.

The function $\log{\phi}:\s\setminus K\fle \R$ is a proper positive function on $\s\setminus K$. Therefore, in order to prove that $\s\setminus K$ is
parabolic it suffices to see that $\log{\phi}$ is superharmonic on $\s\setminus K$. Equivalently, it suffices to see that $\log{\phi}$ is superharmonic on the dense subset
$\s'\subset\s\setminus K$ where it is smooth. In what follows, we will work on that subset $\s'$. From \eqref{normgradh} and \eqref{laph} we get
\beq
\label{laph2}
\Delta h^2=2h\Delta h+2\|\n h\|^2=2(\Theta^2-1).
\eeq
On the other hand, as the function $\bar{r}$ does not depend on $t$ then
$\bar{\n}\bar{r}(x,t)=\hat{\n}\hat{r}(x)$ and $\bar{\Delta}\bar{r}(x,t)=\hat{\Delta}\hat{r}(x)$. Therefore,
\beq
\label{lapr21}
\bar{\Delta}\bar{r}^2(x,t)=2 \bar{r}(x,t)\bar{\Delta}{\bar{r}}(x,t)+2\|\bar{\n}\bar{r}(x,t)\|^2=
2(\hat{r}(x)\hat{\Delta}\hat{r}(x)+1),
\eeq
since, as is well-known, $\|\bar{\n}\bar{r}\|^2 =\|\hat{\n}\hat{r}\|_M ^2=1$.
Applying now Lemma \ref{lemmaaux} to $\psi=r^2$ we get
\beq
\label{lapr22}
\Delta r^2=\bar{\Delta}\bar{r}^2+\hat{\n}^2 \hat{r}^2(N^\ast,N^\ast)=
2(r\hat{\Delta}\hat{r}+1)+\hat{\n}^2 \hat{r}^2(N^\ast,N^\ast).
\eeq

To study the last term of \eqref{lapr22} we will compute first the Hessian of $\hat{r}^2$ on $M$ at a point $x$.
For $v,w \in T_xM$ we have
\begin{eqnarray*}
\hat{\n}^2 \hat{r}^2(v,w) & = & 2\g{\hat{\n}_v (\hat{r}\hat{\n}\hat{r})}{w}_M\\
{} & = & 2\hat{r}(x)\hat{\n}^2\hat{r}(v,w)+2\g{\hat{\n}\hat{r}(x)}{v}_M\g{\hat{\n}\hat{r}(x)}{w}_M.
\end{eqnarray*}
In particular, for $\tau \bot_M \hat{\n}\hat{r}$ of unit length $\|\tau\|_M=1$ we get
\begin{eqnarray*}
\hat{\n}^2\hat{r}^2(\hat{\n}\hat{r}(x),\hat{\n}\hat{r}(x)) & = & 2,\\
\hat{\n}^2 \hat{r}^2(\hat{\n}\hat{r}(x),\tau) & = & 0,\\
\hat{\n}^2 \hat{r}^2(\tau,\tau) & = & 2\hat{r}(x)\hat{\n}^2\hat{r}(\tau,\tau)=2\hat{r}(x)\hat{\Delta}\hat{r}(x).
\end{eqnarray*}
As any $v\in T_xM$ can be decomposed as
\[
v=\g{v}{\hat{\n}\hat{r}(x)}_M\hat{\n}\hat{r}(x)+\g{v}{\tau}_M\tau,
\]
we finally obtain
\[
\hat{\n}^2 \hat{r}^2(v,v)=2\g{v}{\hat{\n}\hat{r}}_M^2+2\hat{r}(x)\hat{\Delta}\hat{r}(x)\g{v}{\tau}_M^2.
\]
Therefore, along the surface $\s'$ we have that
\[
\hat{\n}^2\hat{r}^2(N^\ast,N^\ast)=2\g{N^\ast}{\hat{\n}\hat{r}}_M^2+2r\hat{\Delta}\hat{r}\g{N^\ast}{\tau}_M^2,
\]
and \eqref{lapr22} becomes
\beq
\label{lapr2}
\frac{1}{2}\Delta r^2=r\hat{\Delta}\hat{r}(1+\g{N^\ast}{\tau}_M^2)+1+\g{N^\ast}{\hat{\n}\hat{r}}_M^2.
\eeq

Now, from \eqref{laph2} and \eqref{lapr2} we get that
\begin{eqnarray}
\label{lapphi}
\frac{1}{2}\Delta \phi & = & \frac{1}{2}\Delta r^2-\frac{1}{2}\Delta h^2\\
\nonumber {} & = & r\hat{\Delta}\hat{r}(1+\g{N^\ast}{\tau}_M^2)+\g{N^\ast}{\hat{\n}\hat{r}}_M^2+2-\Theta^2.
\end{eqnarray}
As $M^2$ is complete and has non-negative Gaussian curvature, by the Laplacian comparison theorem we have that
$\hat{\Delta}\hat{r} \leq 1/\hat{r}$, so that
\[
r\hat{\Delta}\hat{r}\leq 1
\]
on $\s'$. Using this in \rf{lapphi}, we obtain that
\begin{equation}
\label{lapphiBIS}
\frac{1}{2}\Delta \phi\leq \|N^\ast\|^2+3-\Theta^2=2,
\end{equation}
since $\|N^\ast\|^2=\Theta^2-1$.

On the other hand, $\nabla\phi=2r\nabla r-2h\nabla h$, and so
\[
\|\n\phi\|^2=4r^2\|\n r\|^2-8rh\g{\n r}{\n h}+4h^2\|\n h\|^2.
\]
Observe that $\bar{\n}\bar{r}=\n r-\g{\bar{\n}\bar{r}}{N}N$ and $\partial_t=-\nabla h-\Theta N$. Taking into account that
$\|\bar{\nabla}\bar{r}\|^2=1$ and $\g{\bar{\nabla}\bar{r}}{\partial_t}=0$, it follows from here that
\[
\|\n r\|^2=1+\g{\bar{\n}\bar{r}}{N}^2, \quad \mathrm{and} \quad
\g{\n r}{\n h}=-\Theta\g{\bar{\n}\bar{r}}{N},
\]
which jointly with \rf{normgradh} implies that
\begin{eqnarray}
\label{gradphi}
\|\n\phi\|^2 & = & 4r^2(1+\g{\bar{\n}\bar{r}}{N}^2)+8rh\Theta\g{\bar{\n}\bar{r}}{N}+4h^2(\Theta^2-1)\\
\nonumber {} & = & 4\phi+4(r\g{\bar{\n}\bar{r}}{N}+h\Theta)^2 \geq 4\phi.
\end{eqnarray}
Therefore, from \eqref{lapphiBIS} and \eqref{gradphi} we finally get
\beq
\label{lap}
\Delta\log\phi=\frac{1}{\phi^2}(\phi\Delta\phi-\|\n\phi\|^2)\leq 0,
\eeq
which means that $\log\phi$ is a proper positive superharmonic function on $\s'$.
Then, $\s'$ is parabolic, and $\s$ is also parabolic.
\end{proof}

It is interesting to look for some natural conditions under which the assumptions of
Theorem \ref{teorlap} are satisfied. In this context, we have the following result.

\begin{proposition}
\label{propwa}
Let $M^2$ be a complete Riemannian surface and let \f\ be a proper spacelike immersion which eventually lies in
\[
\mathcal{W}_a=\{(x,t)\in \N:|t| \leq a \hat{r}(x)\}
\]
for some $0<a<1$. Then the function $\phi=r^2-h^2$ is eventually positive and proper on \s.
\end{proposition}
\begin{corollary}
Let $M^2$ be a complete Riemannian surface with non-negative Gaussian curvature. Then every proper maximal immersion \f\ with non-empty boundary which eventually lies in $\mathcal{W}_a$ for some $0<a<1$, is parabolic.
\end{corollary}
\begin{proof}[Proof of Proposition\rl{propwa}]
Since $f(\s)$ eventually lies in $\mathcal{W}_a$, then there exists a compact set $K \subset \Sigma$ such that
$h^2\leq a^2r^2$ and $\phi=r^2-h^2>a^2r^2-h^2\geq 0$ on $\Sigma \setminus K$. In order to see that $\phi$ is proper, it
suffices to prove that $(\phi|_{\Sigma \setminus K})^{-1}([0,b])$ is compact for every $b>0$. Let
$\bar{\phi}:\mathcal{W}_a\rightarrow\R$ defined by $\bar{\phi}(x,t)=\hat{r}^2(x)-t^2$, so that $\phi |_{\Sigma \setminus K}=\bar{\phi}\circ f |_{\Sigma \setminus K}$. Then $(\phi|_{\Sigma \setminus K})^{-1}([0,b])=f^{-1}(\bar{\phi}^{-1}([0,b]))$. Since $f$ is proper, then it suffices to
prove that $\bar{\phi}^{-1}([0,b])\subset\mathcal{W}_a$ is compact. Observe that for every $(x,t)\in\mathcal{W}_a$ one has
\[
\bar{\phi}(x,t)=\hat{r}^2(x)-t^2\geq\hat{r}^2(x)-a^2\hat{r}^2(x)=(1-a^2)\hat{r}^2(x).
\]
Therefore, for every $(x,t)\in\bar{\phi}^{-1}([0,b])\subset\mathcal{W}_a$ we have that
\[
\hat{r}^2(x)\leq c^2:=\frac{b}{1-a^2}.
\]
This implies that $\bar{\phi}^{-1}([0,b])\subset \bar{B}(x_0,c)\times [-ac,ac]$, where $B(x_0,c)$ denotes the geodesic disc on $M$ of radius $c$ centered at $x_0$ (see Figure 1). Since $\bar{B}(x_0,c)\times [-ac,ac]$ is compact, the result follows.
\begin{figure}[h]
\begin{center}
\includegraphics[width=7cm]{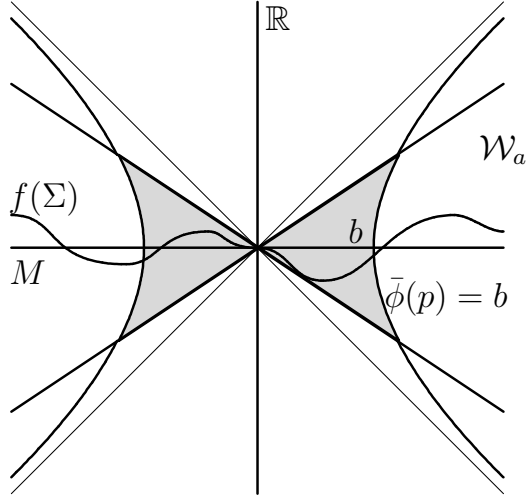}
\caption{$\bar{\phi}^{-1}([0,b])$ is compact under the assumptions of Proposition \ref{propwa}.}
\end{center}
\end{figure}

\end{proof}

On the other hand, recall that a Riemannian surface $(\Sigma,g)$ \textit{without} boundary is
called recurrent if for every nonempty open set $O\subset\Sigma$ with smooth boundary,
$\Sigma\setminus O$ is parabolic (as a surface with boundary). Therefore, as another consequence of our
Theorem\rl{teorlap} we can state the following.
\begin{corollary}
\label{coro5}
Let $M^2$ be a complete Riemannian surface with non-negative Gaussian curvature, and let $\Sigma$ be a maximal surface
in \N without boundary, $\partial\Sigma=\emptyset$. If the function $\phi=r^2-h^2$ is eventually positive and proper
on $\Sigma$, then $\Sigma$ is recurrent.
\end{corollary}
For a proof, simply observe that if $O\subset\Sigma$ is a nonempty open set with smooth boundary, then the function
$\phi$ restricted to $\Sigma\setminus O$ is also eventually positive and proper on $\Sigma\setminus O$, and therefore
the maximal surface with boundary $\Sigma\setminus O$ is parabolic by our main result.

\section{Entire maximal graphs and a Calabi-Bernstein type theorem}

Let $\Omega\subseteq\m^2$ be a connected domain. Every smooth function $u\in\mathcal{C}^\infty(\Omega)$ determines a
graph over $\Omega$ given by $\Sigma(u)=\{ (x,u(x)) : x\in\Omega \}\subset\N$. The metric induced on $\Omega$ from the
Lorentzian metric on the ambient space via $\Sigma(u)$ is given by
\beq
\label{gu}
\g{}{}=\g{}{}_M-du^2.
\eeq
Therefore, $\Sigma(u)$ is a spacelike surface in \N\ if and only if $\|\hat{\n}u\|_M^2<1$ everywhere on $\Omega$.
A graph is said to be entire if $\Omega=\m$.

Consider $\Omega\subseteq\m$ a connected domain and let $x_0\in\mathrm{int}(\Omega)$. We will say that $\Omega$ is
starlike with respect to $x_0$ if for every $x\in\Omega$ there exists a (non-necessarily unique) minimizing geodesic
segment from $x_0$ to $x$ which is contained in $\Omega$. Obviously, if \m\ is a complete Riemannian surface, then
\m\ itself is starlike with respect to any of its points.
\begin{proposition}
\label{lemageom}
Let $M^2$ be a complete Riemannian surface and let $\s(u)$ be a spacelike graph over a domain $\Omega$ which is starlike
with respect to some point $x_0\in\mathrm{int}(\Omega)$.
Then the function $\phi=r^2-h^2$ is eventually positive and proper on $\s(u)$.
\end{proposition}

\begin{proof}
We may assume without loss of generality that $u(x_0)=0$. Since $\s(u)$ is homeomorphic to $\Omega$ (via the standard
embedding $x\in\Omega\hookrightarrow(x,u(x))\in\s(u)$), and the thesis of our result is topological, it is equivalent to prove
that the function $\varphi=\hat{r}^2-u^2$ is eventually positive and proper on $\Omega$. Consider
\[
\mathcal{W}=\{(x,t) \in \Omega \times \R: \hat{r}^2(x)-t^2 \geq 0\}.
\]

Firstly, we will prove that $\varphi$ is positive for every $x\in\Omega-\{x_0\}$. That is, we are going to see that
\[
\label{intW} \s(u)-\{(x_0,0)\} \subset \mathrm{int}(\mathcal{W}).
\]
For a given $x\neq x_0$, consider $\gamma:[0,\ell]\rightarrow\Omega$ a minimizing geodesic segment such that
$\gamma(0)=x_0$, $\gamma(\ell)=x$ and $\ell=\mathrm{dist}_M(x_0,x)=\hat{r}(x)>0$. Let
$\alpha(s)=(\gamma(s),u(s))\in\s(u)$, where $u(s):=u(\gamma(s))$. $\s(u)$ is a spacelike surface, so
$\alpha'(s)=(\gamma'(s),u'(s))\neq (0,0)$ is a non-vanishing spacelike vector, that is,
\[
\g{\alpha'(s)}{\alpha'(s)}=\|\gamma'(s)\|_M^2-u'(s)^2=1-u'(s)^2>0.
\]
Therefore, $-1<u'(s)<1$ for every $0\leq s\leq\ell=\hat{r}(x)$, and integrating we obtain
\beq
\label{u}
-\hat{r}(x)<u(x)<\hat{r}(x).
\eeq
Consequently, $\varphi(x)>0$ and $(x,u(x))\in\mathrm{int}(\mathcal{W})$ for every
$x\in\Omega$, $x\neq x_0$.

It remains to prove that $\varphi$ is proper.
Let us consider on $\m^2\times\R$ the standard Riemannian metric, $\g{}{}_M+dt^2$, and let us denote by
$\mathrm{dist}_+(,)$ the distance related to this Riemannian metric. Let us see now that
\beq
\label{dist1}
\mathrm{dist}_+((x,t),\partial{\mathcal{W}})=
\frac{1}{\sqrt{2}}\min\{\hat{r}(x)-t,\hat{r}(x)+t\}=\frac{1}{\sqrt{2}}(\hat{r}(x)-|t|)
\eeq
for every $(x,t)\in\mathcal{W}$. Observe that $\partial{\mathcal{W}}$ can be decomposed into
$\partial{\mathcal{W}}=\partial{\mathcal{W}^+} \cup \partial{\mathcal{W}^-}$ where
\[
\partial{\mathcal{W}^+}=\{(x,\hat{r}(x)):x \in \Omega\} \quad
\mathrm{and} \quad \partial{\mathcal{W}^-}=\{(x,-\hat{r}(x)):x \in \Omega\}.
\]
Therefore,
\[
\mathrm{dist}_+((x,t),\partial{\mathcal{W}})=
\min\{\mathrm{dist}_+((x,t),\partial{\mathcal{W}^+}),\mathrm{dist}_+((x,t),\partial{\mathcal{W}^-})\}.
\]
Expression \rf{dist1} is clear for $x=x_0$ (and necessarily $t=0$). For a given $x\neq x_0$, let
$\gamma:[0,\hat{r}(x)]\rightarrow\Omega$ be a minimizing geodesic segment such that $\gamma(0)=x_0$,
$\gamma(\hat{r}(x))=x$. We will compute first $\mathrm{dist}_+((x,t),\partial{\mathcal{W}^+})$. Since $\gamma$ is
minimizing, for every $s\in[0,\hat{r}(x)]$ we have $\hat{r}(\gamma(s))=s$, so that
$(\gamma(s),s)\in\partial{\mathcal{W}^+}$ and
\[
\mathrm{dist}_+((x,t),(\gamma(s),s))^2=
\mathrm{dist}_M(x,\gamma(s))^2+(t-s)^2=
(\hat{r}(x)-s)^2+(t-s)^2.
\]
Observe that this expression attains its minimum at $s_0=(\hat{r}(x)+t)/2$, and
\beq
\label{dist3}
\mathrm{dist}_+((x,t),(\gamma(s_0),s_0))=\frac{1}{\sqrt{2}}|\hat{r}(x)-t|=
\frac{1}{\sqrt{2}}(\hat{r}(x)-t).
\eeq
We claim that $\mathrm{dist}_+((x,t),\partial{\mathcal{W}^+})$ is given by \rf{dist3}. In fact, for every
$y\in\Omega$ we have that
\[
\mathrm{dist}_M(x,y)\geq|\mathrm{dist}_M(x_0,x)-\mathrm{dist}_M(x_0,y)|=|\hat{r}(x)-\hat{r}(y)|,
\]
and so
\begin{eqnarray*}
\mathrm{dist}_+((x,t),(y,\hat{r}(y)))^2 & = & \mathrm{dist}_M(x,y)^2+(t-\hat{r}(y))^2\\
{} & \geq & (\hat{r}(x)-\hat{r}(y))^2+(t-\hat{r}(y))^2\\
{} & \geq & \min_{s\geq 0}\left((\hat{r}(x)-s)^2+(t-s)^2\right)=\frac{1}{2}(\hat{r}(x)-t)^2.
\end{eqnarray*}
Therefore,
\beq
\label{dist4}
\mathrm{dist}_+((x,t),\partial{\mathcal{W}^+})=\frac{1}{\sqrt{2}}(\hat{r}(x)-t).
\eeq

With a similar argument for $\partial{\mathcal{W}^-}$,
\beq
\label{dist5}
\mathrm{dist}_+((x,t),\partial{\mathcal{W}^-})=\frac{1}{\sqrt{2}}(\hat{r}(x)+t).
\eeq
Thus, \eqref{dist1} follows from \eqref{dist4} and \eqref{dist5}.

Let $x\in\Omega$, $x\neq x_0$, and let $\gamma:[0,\hat{r}(x)]\rightarrow\Omega$ be a
minimizing geodesic segment such that $\gamma(0)=x_0$, $\gamma(\hat{r}(x))=x$. Write $u(s)=u(\gamma(s))$. Then
$(\gamma(s),u(s))\in\mathcal{W}$ and by \rf{dist1} we have that
\beq
\label{dist2}
\mathrm{dist}_+((\gamma(s),u(s)),\partial{\mathcal{W}})=
\frac{1}{\sqrt{2}}(s-|u(s)|).
\eeq
This implies that $\mathrm{dist}_+((\gamma(s),u(s)),\partial{\mathcal{W}})$ is a positive increasing function for
$0<s\leq\hat{r}(x)$. Therefore, if we choose $\delta>0$ such that the geodesic disc
$B_\delta=B(x_0,\delta)\subset\subset\Omega$, then it follows that
\beq
\label{dist6}
\mathrm{dist}_+((x,u(x)),\partial{\mathcal{W}})\geq\varepsilon>0 \quad \textrm{ for every } x\in\Omega\setminus B_\delta,
\eeq
where
\[
\varepsilon=\min_{x\in\partial B_\delta}
\mathrm{dist}_+((x,u(x)),\partial{\mathcal{W}}) >0
\]
(see Figure 2).

\begin{figure}[h]
\begin{center}
\includegraphics[width=7cm]{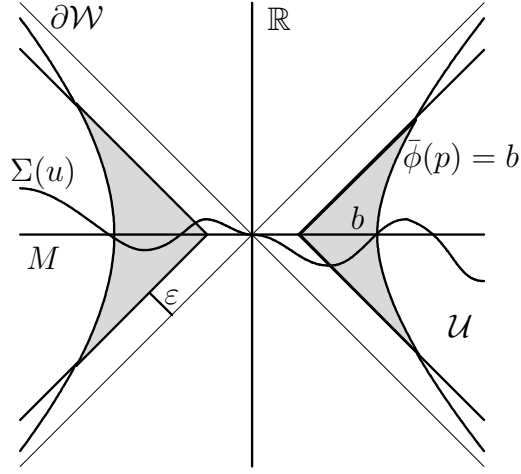}
\caption{$\mathrm{dist}_+((x,t),\partial \mathcal{W}) \geq \varepsilon$ eventually in $\Omega$.}
\end{center}
\end{figure}

Now we are ready to prove that $\varphi$ is proper on $\Omega$. Since $\Omega=\overline{B_\delta}\cup(\Omega\setminus B_\delta)$ with
$\overline{B_\delta}$ compact, it suffices to prove that $\varphi|_{\Omega\setminus B_\delta}$ is
proper on $\Omega\setminus B_\delta$. Let $f:\Omega\rightarrow\Omega\times\R$ be the standard
embedding, $f(x)=(x,u(x))$, and let $\bar{\phi}:\Omega\times\R\rightarrow\R$ defined by $\bar{\phi}(x,t)=\hat{r}^2(x)-t^2$.
Observe that $f$ is trivially proper for if $A\subset\Omega\times\R$ is compact, then
\[
f^{-1}(A)\subset f^{-1}(\pi_1(A)\times\pi_2(A))=\pi_1(A)\cap u^{-1}(\pi_2(A))
\]
is compact, where $\pi_1:\Omega\times\R\rightarrow\Omega$ and $\pi_2:\Omega\times\R\rightarrow\R$ are the projections.
By \rf{dist6} we have that
\[
f(\Omega\setminus B_\delta)\subset\mathcal{U}=\{ (x,t)\in\mathcal{W} :
\mathrm{dist}_+((x,t),\partial{\mathcal{W}})\geq\varepsilon \}.
\]
Therefore,
$\varphi|_{\Omega\setminus B_\delta}=\bar{\phi}|_{\mathcal{U}}\circ f|_{\Omega\setminus B_\delta}$. Note
that the map $f|_{\Omega\setminus B_\delta}:\Omega\setminus B_\delta\rightarrow\mathcal{U}$ is proper. In fact, for every
compact set $A\subset\mathcal{U}$ we see that
$(f|_{\Omega\setminus B_\delta})^{-1}(A)=f^{-1}(A)\cap(\Omega\setminus B_\delta)$ is compact. Therefore, it suffices to
show that $\bar{\phi}|_{\mathcal{U}}:\mathcal{U}\rightarrow\R$ is proper, or, equivalently, that for every $b>0$,
$(\bar{\phi}|_{\mathcal{U}})^{-1}([0,b])=\mathcal{U}\cap\bar{\phi}^{-1}([0,b])$ is compact. Let
$(x,t)\in\mathcal{U}\cap\bar{\phi}^{-1}([0,b])$. Since $(x,t)\in\mathcal{U}$, by \rf{dist1} we obtain that
\[
|t|\leq\hat{r}(x)-\sqrt{2}\varepsilon.
\]
Therefore, since $\bar{\phi}(x,t)\leq b$, this gives
\[
\hat{r}^2(x)-b\leq t^2\leq(\hat{r}(x)-\sqrt{2}\varepsilon)^2,
\]
that is
\[
\hat{r}(x)\leq c:=\frac{2\varepsilon^2+b}{2\sqrt{2}\varepsilon}.
\]
This implies that
$(\bar{\phi}|_{\mathcal{U}})^{-1}([0,b])\subset\overline{B}(x_0,c)\times[\sqrt{2}\varepsilon-c,c-\sqrt{2}\varepsilon]$
is compact. This finishes the proof of Proposition\rl{lemageom}.
\end{proof}

As another consequence of our main result (Theorem\rl{teorlap}), we can give a new proof of the following Calabi-Bernstein
theorem, first established in \cite[Theorem 4.3]{AA} (see also \cite{AA2} for another approach to the parametric version of
that result, first established in \cite[Theorem 3.3]{AA}).
\begin{corollary}
\label{CB}
Let $\m^2$ be a complete Riemannian surface with non-negative Gaussian curvature. Then any
entire maximal graph $\Sigma(u)$ in \N\ is totally geodesic. In addition, if $K_M>0$ at some point on \m, then $u$ is
constant.
\end{corollary}
\begin{proof}
By applying Proposition\rl{lemageom} with $\Omega=\m$ we know that the function $\phi=r^2-h^2$ is eventually
positive and proper on $\Sigma(u)$. Therefore, by Corollary\rl{coro5} we have that $\Sigma(u)$ is recurrent. Equivalently,
any negative subharmonic function on the surface $\Sigma(u)$ must be constant. The proof then follows as in the proof
of \cite[Theorem 3.3]{AA}. For the sake of completeness, we briefly sketch it here. Recall that
$\Theta=\g{N}{\partial_t}\leq -1<0$. A careful computation using Codazzi equation and the maximality of $\Sigma(u)$ gives
that
$$
\Delta\Theta=\Theta(\kappa_M(\Theta^2-1)+\|A\|^2)
$$
and
$$
\|\nabla\Theta\|^2=\frac{1}{2}\|A\|^2\|\nabla h\|^2=\frac{1}{2}\|A\|^2(\Theta^2-1),
$$
where $\|A\|^2=\mathrm{tr}(A^2)$ and $\kappa_M$ stands for the Gaussian curvature of \m\ along $\Sigma(u)$.
This implies that
$$
\Delta\left(\frac{1}{\Theta}\right)=-\frac{\Delta\Theta}{\Theta^2}+\frac{2\|\nabla\Theta\|^2}{\Theta^3}=
-\frac{1}{\Theta}\left(\kappa_M(\Theta^2-1)+\frac{\|A\|^2}{\Theta^2}\right)\geq 0.
$$
That is, $1/\Theta$ is a negative subharmonic function on the recurrent surface $\Sigma(u)$, and hence it must be
constant. Therefore, $\|A\|^2=0$ and $\kappa_M(\Theta_0^2-1)=0$ on $\Sigma(u)$, where $\Theta=\Theta_0\leq -1$.
Thus, $\Sigma(u)$ is totally geodesic and, if $\kappa_M>0$ at some point on $\Sigma(u)$ (equivalently, $K_M>0$ at some point on
\m), then it must be $\Theta_0=-1$, which means that $h$ is constant and $\Sigma(u)$ is a slice.

\end{proof}

\bibliographystyle{amsplain}

\end{document}